# A Conversation with Yuan Shih Chow

**Zhiliang Ying and Cun-Hui Zhang**

*Abstract.* Yuan Shih Chow was born in Hubei province in China, on September 1, 1924. The eldest child of a local militia and political leader, he grew up in war and turmoil. His hometown was on the front line during most of the Japanese invasion and occupation of China. When he was 16, Y. S. Chow journeyed, mostly on foot, to Chongqing (Chung-King), the wartime Chinese capital, to finish his high school education. When the Communist party gained power in China, Y. S. Chow had already followed his university job to Taiwan. In Taiwan, he taught mathematics as an assistant at National Taiwan University until he came to the United States in 1954. At the University of Illinois, he studied under J. L. Doob and received his Ph.D. in 1958. He served as a staff mathematician and adjunct faculty at the IBM Watson Research Laboratory and Columbia University from 1959 to 1962. He was a member of the Statistics Department at Purdue University from 1962 to 1968. From 1968 until his retirement in 1993, Y. S. Chow served as Professor of Mathematical Statistics at Columbia University. At different times, he was a visiting professor at the University of California at Berkeley, University of Heidelberg (Germany) and the National Central University, Taiwan. He served as Director of the Institute of Mathematics of Academia Sinica, Taiwan, and Director of the Center of Applied Statistics at Nankai University, Tianjin, China. He was instrumental in establishing the Institute of Statistics of Academia Sinica in Taiwan. He is currently Professor Emeritus at Columbia University.

Y. S. Chow is a fellow of the Institute of Mathematical Statistics, a member of the International Statistical Institute and a member of Taiwan's Academia Sinica. He has numerous publications, including *Great Expectations*: *The Theory of Optimal Stopping* (1971), in collaboration with Herbert Robbins and David Siegmund, and *Probability Theory* (1978), in collaboration with Henry Teicher. Y. S. Chow has a strong interest in mathematics education. He taught high school mathematics for one year in 1947 and wrote a book on high school algebra in collaboration with J. H. Teng and M. L. Chu. In 1992, Y. S. Chow, together with I. S. Chang and W. C. Ho, established the Chinese Institute of Probability and Statistics in Taiwan.

This conversation took place in the fall of 2003 in Dobbs Ferry, New York.

*Zhiliang Ying is Professor and Co-Chair, Department of Statistics, Columbia University, New York, New York 10027, USA. Cun-Hui Zhang is Professor of Statistics, Rutgers University, New Brunswick, New Jersey 08854-8019, USA.*

## FAMILY BACKGROUND AND EARLY YEARS

**Query:** We have always been fascinated by the interesting stories you have told us about yourself, but there are certainly a lot more. It would be greatly appreciated if you could share your experience and thoughts with colleagues, students, and friends in and outside the statistical community. We wonder if you could start with your childhood.

**Chow:** I was born in 1924 into a family with a large estate in a rural town in central China. It was a time of war and turmoil. My family had been dependent on land for our livelihood since the Qing Dynasty. My father, a local militia leader, had more than 100 militiamen to protect our hometown. When I was four, my father lost a battle against the communists and our home was burned down. We first left for Wuhan, but my father soon found the provincial capital too expensive for an extended family of 40 people without income from land. So we eventually settled in Xiangfan, a city near my hometown.

**Query:** How did you become interested in mathematics?

**Chow:** During the period of our "refugee" life, the family hired an old-fashioned tutor, who taught traditional Chinese language and literature. I started formal schooling in a highly scholastic primary school when I was eight, after we settled down in Xiangfan. I started as a second grader, but I was not able to recognize Arabic numerals. I got a big shock when I received a "0" mark on my first math exam, but it motivated me to work hard. My efforts paid off, resulting in excellent grades in subsequent exams. That is how I got interested in mathematics.

## MIDDLE SCHOOLS AND WORLD WAR II

**Query:** Your teenage years overlapped with the Sino–Japanese War. What happened to your education?

**Chow:** The war broke out in 1937, right after my graduation from the primary school. Soon schools were all shut down when the Japanese began bombing our area. Although I passed the entrance examination for a provincial middle school, I had to wait until the end of 1938 to attend a new combined junior high school.

**Query:** Were you still able to get a decent education?

**Chow:** The combined junior high was a good school, but the combined senior high lacked qualified teachers since many had left for mountainous areas to avoid the war. Upon completion of the junior high, three classmates and I decided to find a better high school outside our area. That turned out to be quite an adventure.

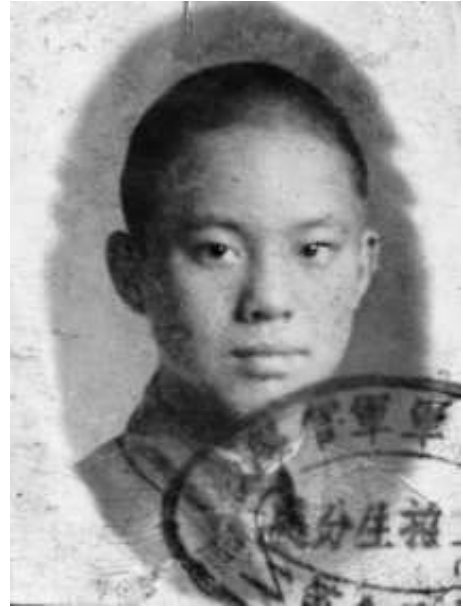

Fig. 1.  *Y. S. Chow in 1940.*

We first went to the wartime provincial capital. We walked for two whole weeks through mountains, including the Three Gorges, in the midst of war and cold winter weather. It was quite a dangerous trip. However, to our disappointment we were "drafted" into an agricultural school there, which was treated as a military school during the war.

**Query:** What happened then? Did you stay there?

**Chow:** Because it was a "military school," running away was punishable by imprisonment. Despite this, in July 1941, we decided to run away for Chongqing, the wartime capital of China. We left the school at night and walked for four days to arrive at a port city on the Yangtze River. By then we were penniless. We spent a night in a park under the open sky, and traded a gold ring for our boat passage to Chongqing the next morning. In Chongqing, we stayed in a youth center run by the government, which also served as a shelter for refugee students.

**Query:** What did you do in Chongqing?

**Chow:** We decided to take an entrance exam for the No. 2 National Middle School, located about 90 kilometers, or one day and a half by walking, outside of Chongqing. From the youth center, we got a small amount of money for the cost of the free meals we would miss for three days during our journey. However, upon our arrival, we learned that we had to arrive one day earlier to apply. We pleaded with the principal and argued that because we could afford food for only three days, we really were not able to come earlier. The principal finally allowed us to take the exam. He even let us sleep on the tables in the dining hall. I passed the exam, which was quite an honor and resulted in a cash award from the school district of my hometown. Better yet, I did not have to walk back to Chongqing!

**Query:** Were you able to get a good education there, especially in mathematics?

**Chow:** The No. 2 National Middle School, a combined school of some of the best schools around Shanghai, was one of the most envied among all the key schools in China during the war. The students there were all refugees from faraway, but they enjoyed very high status and were very proud of themselves. We had excellent teachers. Among the best was Mr. Mao, our mathematics instructor. Mr. Mao had a great influence on the way we studied mathematics. He emphasized careful reading of advanced textbooks, not working on numerous exercises. Students in his class enjoyed total freedom: free to decide whether to do homework and whether to attend classes! However, they all worked hard and excelled in their studies. The percentage of graduates admitted into top universities was very high. Mr. Mao was also the student advisor of our class. On holidays, he would invite us to his house for dinner, since we had no place else to go.

## ZHEJIANG UNIVERSITY

**Query:** It is impossible to imagine nowadays what you went through seeking a good education at such young ages. What happened during your college years?

**Chow:** During the war, many universities moved with the central government to southwest China. Since our school was the best near Chongqing, they often sent professors there to recruit students. A common practice was to admit top students without entrance examinations. Zhejiang University was the first to come during my senior year. When I received an offer from it, I decided to seize the opportunity.

**Query:** Why did you choose mathematics as your major?

**Chow:** I had a childhood dream that someday I would found a private middle school in my own hometown. My family had approximately 1000 *mus* (about 200 acres) of land. I could draw out a certain portion of income from the land to keep the school going. Students could work part time on the land so they would not have to pay tuition. The main problem would be to hire qualified teachers, especially mathematics teachers. This motivated me to choose mathematics as my major.

**Query:** What was your college life like during the war?

**Chow:** In the summer of 1944, I left Chongqing by bus for the wartime campus of Zhejiang University. The trip was dangerous, as accidents were common on the bus route through the mountain paths. The bus, which burned charcoal to generate steam, was in terrible condition. Shortly after our arrival, in December 1944, the Japanese reached areas near the campus. School had to be suspended. Since we were worried that the Japanese could interrupt our financial support from the central government, our first responsibility was to find food. The second was to protect female students. Each student got a gun. Meanwhile, we organized Delegates to Battle Fields to visit the front to prop up the morale of soldiers. A few months later, the Yunnan–Burma road opened. We saw soldiers in large U.S. made ten-wheel trucks heading toward central China. We all cheered, "Americans have come to help! There is hope at last!" Although we did not see any American soldiers, that was a turning point in the Sino–Japanese War. I still vividly remember the end of the war. It was a hot summer day. We were swimming and heard someone crying, "Japanese have surrendered!" We were overjoyed. That was the beginning of my second year at Zhejiang University. A year later, the university moved back to its permanent campus in Hangzhou, the capital of Zhejiang province.

**Query:** Who were your teachers at Zhejiang University? How did they influence your life?

**Chow:** The most senior faculty members were Buqing Su and Jiangong Chen. When we returned to Hangzhou, Chen took a position in Academia Sinica in Shanghai. Ruiyun Xu, a student of Carathéodory, directed my thesis. She wrote an excellent book on real variables and later became Chairwoman of the Department of Mathematics at Hangzhou University.

The Chinese Civil War broke out shortly after World War II. When I was a junior, my hometown changed hands repeatedly between the Communist troops and the Nationalist government, so my family could no longer support me financially. I sought help from Su, who used his connection to secure a teaching job for me at a middle school. Su was very kind to students, but he was also a very strict teacher. Su's lectures were very serious. He would not speak a word before all the students arrived, so nobody dared to be tardy. Chen and Su were the most respected professors by colleagues and students alike at Zhejiang University. Every year, Chen and Su gave a formal party in a fine restaurant for students. Everyone had to have a drink to the health of each teacher. Since we had about a dozen teachers, we all learned to drink.

## TAIWAN

**Query:** What happened after you graduated from college?

**Chow:** In 1947 and 1948, Su and Chen went to Taiwan for the formation of the National Taiwan University, which was Taipei Imperial University under Japanese rule. Chen became the first Provost of National Taiwan University and Su became the first Dean of the College of Science. Su also founded the mathematics department there. Because I took a one-semester leave to teach in the middle school, I graduated in January 1949 instead of July 1948. I planned to do my postgraduate work with Chen. The department waived my entrance examination, but the graduate program had only fall enrollment. What was I supposed to do to fill up the gap? Chen suggested that I spend six months in Taiwan to help the mathematics department there. However, I did not expect that I would not be able to return to the Mainland. Soon after I left for Taiwan, the Nationalist government was defeated and retreated to Taiwan.

I still dreamed about pursuing postgraduate study. At that time, students in Taiwan were not allowed to go overseas, especially male students since they had to be ready for military service if called upon. Later, the government began to lift the restrictions. H. C. Chow, a mathematician, suggested to the government that Taiwan should adopt an open policy so that students would have hope about their future; otherwise they could become a source of political instability. Specifically, he suggested that the government should allow a limited number of students to go abroad after passing certain examinations. The government accepted the suggestion. As a result, students became more motivated in their study. I believe this indeed helped the government to prevent student movements. The first exam took place in 1953, but I had no connection to learn about the exam in advance, so I came to America in 1954.

The Mathematics Department at Zhejiang University gave very low grades. Even though I did well, my transcript was not presentable. Fortunately, Chung-Tao Yang, who was two years ahead of me at Zhejiang, was then a postdoc at Illinois. He wrote a letter for me saying that he and I got exactly the same training at Zhejiang and that I would have no problem completing the Ph.D. program at Illinois. With his reference, I got a fellowship award from Illinois. I learned from this experience so that I was always generous in marking the papers of my students.

## URBANA–CHAMPAIGN

**Query:** So, you came to Illinois. Why did you decide to study probability with Doob?

**Chow:** I arrived at Champaign more than two weeks late. It took me 33 days by ship to reach Oregon, and then a three-day ride on a Greyhound bus. The steamer ticket, costing 100 U.S. dollars after a 50% discount for students, was not so easy to purchase. Red tape and mandatory political training delayed my passport application.

I first went to see H. Ray Brahana, Chairman of the Mathematics Department. After examining my transcript, he found that I had taken differential geometry, real variables and other advanced courses as an undergraduate. He concluded that I had an extremely solid foundation and put me in all the hard courses least selected by graduate students: topology, Riemann geometry, measure theory and abstract algebra (taught by himself). At first, I had a hard time coping with the course materials, due in part to my problem with English. I could not remember how I survived, but my grades for midterm exams turned out fine. I felt better and began to blend in with my classmates. I made many good friends, including Jim Abbott.

At the end of my first year, there was a summer camp in Indiana for graduate students from China. The camp invited Shiing-Shen Chern, then at the University of Chicago, to give a lecture. After his talk, I approached him for advice about choosing a specialty and a Ph.D. advisor. Chern told me to go with Doob if he was willing to accept me. When I returned to Illinois, the first thing I did was to find Doob, whom I had never met and knew very little about. Doob said O.K. and asked me what I intended to work on. I replied Fourier series, since I had studied the topic as an undergraduate. "Then, you should be my advisor," he joked. When I asked him what I should learn from him, he answered, "Probability."

**Query:** What was it like with your thesis research?

**Chow:** I took Doob's probability course in 1955, after I chose him as my advisor. There were four students in the class, including Jim Abbott and me. It was actually one of the largest probability classes. Among the auditors were Don Burkholder, who had just arrived at Illinois as an Assistant Professor, and Klaus Krickeberg from Germany.

From 1956 to 1957, Doob took a sabbatical leave in Europe. I wrote a paper on martingales while he was away. I felt very satisfied. After I mailed him the paper, I let myself relax and really had a good time. However, one day I realized that I still had to prove a set-theoretical lemma upon which my paper was based. I could not prove it although I had thought that it was obvious when I first wrote it down. I had to let Doob know the situation. He asked me if there was a counterexample. I did not have one, but I just could not prove the lemma. I decided to start afresh with a different approach. Krickeberg was also working on the topic. His advisor C. Y. Pauc was visiting Purdue University. Since Doob was in Europe, I went to see Pauc often. In this sense, Pauc was my quasi-advisor. I finished my thesis work when Doob came back. Two months later, in January 1958, I got my Ph.D. degree.

**Query:** Besides probability, were you involved in any other research project?

**Chow:** Illinois required a minor, so I selected physics. It was an easy decision. I had taken Riemann geometry in my first year, which required relativity theory. To fulfill the physics minor, I just needed to take relativity and hydrodynamics, which were very close. My advisor for my minor was A. H. Taub, who worked on the first computer at Princeton with von Neumann.

After sending my paper to Doob in Europe, I sought help from Taub to find a job for me. He gave me one right away. It was a Research Associate position with his NASA project. Sputnik had just gone into the space. Taub said it was one of the most important and urgent projects in the country. He offered to help me apply for a green card, but he also demanded that I stop working on probability if I took the job. It was very hard then for a Chinese person to get a green card, since the quota for China was only 105 per year, compared with 20,000 for each European nation. The job offer was hard to resist. I could get a green card while working on an important problem and I liked Illinois. When Doob came back, I handed him the thesis and told him that I would have to say goodbye to probability, not even to publish my thesis. I jumped into Taub's computer lab.

The specific problem I worked on came from hydrodynamics. The problem has five differential equations and six unknown quantities. By 1959, after over a year's effort, I came up with a principle, which could lead to a sixth equation. Taub was delighted. He said he would present the result in a talk in Europe and asked me to compute the solution. The capacity of our computer was 1K, so that we could compute the solution at no more than 50 points. It turned out that the domain of our solution was very narrow, requiring at least 250 points. We had to wait for a more powerful computer!

## IBM AND COLUMBIA

**Query:** I remember you told me that you would probably also do very well if you stayed with computers. What else did you do in that field?

**Chow:** One day in 1959, IBM sent a recruiter to Illinois. He told me that I was exactly the man they were looking for. IBM had just begun to form the Watson Research Center, and it had appointed Herman Goldstine as Director of Scientific Development in its Data Processing Division. Goldstine was a deputy director of the computer lab under von Neumann in the Institute for Advanced Studies. After von Neumann's death, the Institute dissolved its computer group. As the head of a new division, Goldstine certainly wanted to hire. I was very knowledgeable with software and computer operations, and I was a student of Taub's, so I was a perfect match for IBM. I was pleased to get IBM's offer.

Since I was leaving Illinois anyway, I decided to explore the job market. Jacob Wolfowitz made offers to both Roger Farrell and me. Farrell, my contemporary at University of Illinois, accepted the offer, but I decided to choose IBM. I liked New York City. Cornell was an excellent university, but Ithaca seemed rather isolated. It did not hurt that IBM's salary was double what Cornell offered.

Since IBM had not yet built the Yorktown campus for the Watson Center, it put Goldstine's division at the Lamb estate near New York City. At that time, Purdue was the only university with a department of mathematics and statistics, so Goldstine lured Purdue's chair, Carl Kossack, to be in charge of the statistics group. When I arrived at IBM, Kossack asked me if I wanted to join his group, but I said statistics was not my field. Instead, I joined the pattern recognition group, which included two other people: one from Princeton specializing in topology and the other from Harvard specializing in physics. William Feller was our consultant. Our main assignment was to study the feasibility of developing computer-based methodology to read the Soviet newspaper *Pravda* and translate its content into English. That is, we wished to be able to place a copy of *Pravda* in the computer and automatically produce its English version. We spent a lot of time and eventually decided that the project would take at least 10 years. Because of our report, IBM did not participate in the project, while other big corporations, including GE and RCA, did and lost millions of dollars. This might be my greatest contribution to IBM.

**Query:** How did you return to probability?

**Chow:** After we submitted our report, I became idle again and began to think about my future. I had promised Taub not to work on probability. I had worked in hydrodynamics with him and then pattern recognition with IBM. However, I was no longer working for Taub and thus free to return to probability. I decided to attend seminars at Columbia University.

Ron Pyke was the seminar chair. After the first seminar, I introduced myself and asked Ron if I could attend the seminars. He said I was very welcome. Herb Robbins heard our conversation and asked me if I had to go back to IBM right after the seminar. When I replied no, he took me to his office and said that it was too bad I missed his more interesting seminar a week earlier. He went on to tell me that it was a secretary problem and provided details on the blackboard.

I worked on Robbins' problem afterward, using martingales. One week later, I went to see Robbins to show him my result. He went over it and asked me what he was supposed to do with it. When I realized he was talking about a joint paper, I told him that it was his stuff. However, Robbins said he did not prove his theorem and my proof was more general. This was the beginning of our collaboration. It was around the end of 1959.

**Query:** So how often did you go to Columbia?

**Chow:** I went to Columbia every week after Robbins and I started to collaborate. IBM had a unit called Watson Scientific Computing Laboratory at Columbia University, located in a small house right at 116th Street between Broadway and Claremont Avenue, with a staff of about a dozen people. Wallace Eckert, founder and Director of the Watson Lab, was a very good friend of Robbins'. Robbins felt it was too much trouble for me to come back and forth so often, so he called Eckert to suggest that IBM transfer me over to Manhattan. At that time, Columbia did not have its own computer. IBM provided computing services to Columbia through the Watson Lab. IBM provided the staff, while Columbia took care of the maintenance. The IBM staff members usually had adjunct titles at Columbia. In fact, they were in charge of the applied math program at Columbia. Peter Welch got his Ph.D. degree in statistics while being a student staff member at the Watson Lab.

At the Watson Lab, my duty was to offer a course in statistics at Columbia, and that was all. IBM explained to me that the Watson Lab was a window for IBM, through which IBM could keep informed of what was going on in the outside world. IBM later merged the Watson Lab into Watson Research Center, which was much larger in scale. The Watson Lab had many high-caliber people. It covered many fields, including chemistry, mathematics, physics and astrophysics.

**Query:** Who else was at Columbia at that time? What did you teach?

**Chow:** In the spring semester of 1960, I taught as an Adjunct Assistant Professor at Columbia. Ron Pyke, Don Ylvisaker, Jerry Sacks and Lajos Takacs were Assistant Professors. Joe Gani was probably a Visiting Associate Professor. These were the young people there. I suggested that the department should offer probability, but we needed to offer measure theory first. Ron agreed, so it was the first course I taught at Columbia. We called it "mathematics for statistics," since measure theory was in the territory of the mathematics department. The second course I offered was martingales, which attracted more people. Robbins sat in, as did S. Moriguti and G. Mariyama; both were visiting from Tokyo University.

**Query:** You started teaching at Columbia and working with Robbins before you went to Purdue. Why did you decide to go there?

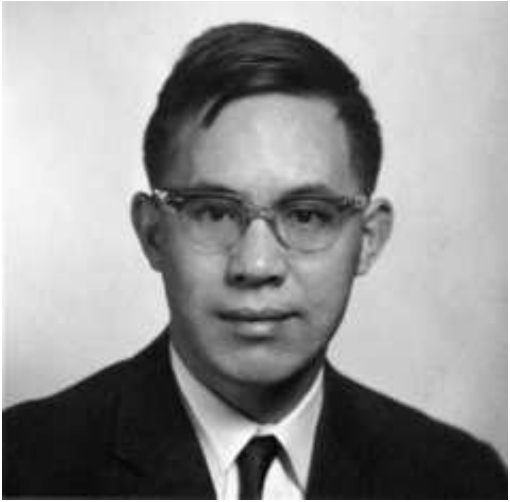

Fig. 2. *Y. S. Chow in 1962.*

**Chow:** In 1961, I was promoted to Adjunct Associate Professor, but IBM paid my salary anyway. IBM staff had to go through a formal evaluation process periodically. When Eckert asked Robbins about his impression of me, Robbins commented that I was more suitable for Columbia than for IBM. Eckert summoned me to his office and told me about his conversation with Robbins, but he said he could not figure out what Robbins actually intended to tell him.

In the same year, I went to Washington, DC, to attend an IMS meeting. I invited Joe Doob and Frank Spitzer for dinner. During our conversation, I mentioned to them that I was working for IBM, but my duty was just to teach a class at Columbia. Spitzer said that it was "too good to last," and his remark kept on bothering me afterward. When I asked Robbins what he meant about his conversation with Eckert, he said he was ready to offer me a job. The idea seemed fine. I would become an Associate Professor automatically. However, the Statistics Department at Columbia did not have a tenure-track line, so I was appointed as a Visiting Associate Professor. Robbins promised that my position would last indefinitely as long as Columbia had money. That was 1962.

After I took the position at Columbia, Shiing-Shen Chern happened to pass by New York. I had kept in touch with him since he advised me to work with Doob. He asked me how I had been doing. I said that I had left IBM to teach at Columbia. He remarked that the job at IBM was a very good one, better than a visiting position with year-to-year appointment. When I told him Robbins' promise, Chern commented, "What difference does that make?" It made me uneasy again.

Henry Teicher, who had been at Purdue since my student days at Illinois, was visiting New York University at that time. I knew him very well because Illinois and Purdue had a get-together every year. When I was a student at Illinois, he was an Associate Professor at Purdue. Since I was Doob's student, it was very easy for me to get to know people well. Henry had graduated from Columbia, started with Stanford and then went on to Purdue. When he visited NYU, he came to Columbia often to hang around with us and I played Go with him. I told him that I had resigned from IBM and was currently looking for a job, since my position at Columbia was not a permanent one. He called Felix Haas at Purdue, who immediately made me an offer: Associate Professor with tenure! I was very happy to get an offer without even having to apply, so I accepted. It was almost too good to be true. Since I was to visit Berkeley the following year with Robbins during Jerzy Neyman and Michel Loève's sabbatical leave, my appointment at Purdue would start in the fall of 1963.

**Query:** So you went to Berkeley with Robbins. What did you do there?

**Chow:** Robbins and I stayed in Neyman's house and shared the same office. Although we should have been able to see each other all the time, it was just the opposite. When I got up early in the morning, he would be fast asleep after coming home way beyond midnight. Of course, when I came back from the office, he was never in the house.

We were working together on the same project. The problem was optimal stopping for the relative frequency with coin tossing. One day, David Blackwell asked me what I was doing. I told him the problem and said we had no idea how to proceed. He said that it was very interesting and he would work on nothing else but this problem. I told Robbins about Blackwell's competing effort and suggested that we work harder to solve the problem first. However, during our entire visit to Berkeley, Robbins spent a lot of time dating girls. Blackwell was a serious researcher. He would come to our office, claiming the existence of the optimal rule at one moment and trying to give a counterexample the next. We were confused. By April or May, I found a proof for the existence of the optimal stopping rule. When Blackwell heard this, he focused his attention on finding the stopping boundary explicitly. It turned out that the square-root boundary is optimal.

I like to have conversations about mathematics and to collaborate with others. In a large sense, the purpose of research is to learn and to understand. Solving problems is far more important than who solves them. This attitude makes it easy for me to work with others.

## PURDUE AND TWO BOOKS

**Query:** So, you went to Purdue after visiting Berkeley. Could you share with us your experience there?

**Chow:** Felix Haas was Chairman of the Department of Mathematics and Statistics at Purdue. But soon after I accepted his offer, Purdue decided to make the department even larger by establishing a new Division of Mathematical Sciences, including the Departments of Mathematics, Statistics and Computer Science. Shanti Gupta was hired from Bell Labs as head of the Statistics Department. When I arrived at Purdue, Gupta was also new there. Teicher had already become a full professor.

I had a joint appointment, half mathematics and half statistics. Traditionally, measure theory was a prerequisite for probability theory. Since measure theory was very hard, most students would be able to start to work on a thesis in functional analysis when they took it. That was why Doob did not have many students at Illinois. I proposed to offer probability theory without requiring measure theory. My colleagues in the Mathematics Department were not sure about my approach. I said I would write my own lecture notes. When they asked why, I replied that I would never get any students otherwise. This is how I started to teach probability without requiring measure theory. My first Ph.D. student was William Stout, whose undergraduate major was electrical engineering. His successful career demonstrates that my approach worked.

I offered my new course and started to write notes. A well-known publisher heard about this and proposed to publish a book based on my notes. However, writing a book is a different matter. My officemate at Purdue was Teicher, so I asked him if he would help by being a co-author. He agreed. The publisher gave us a very attractive offer with a $7000 advance. We took the money without having to write down a single word, and I spent it on a new car. We were supposed to finish the book in three years. However, I left Purdue to visit Columbia in 1965. One year later, Teicher also left, first for Columbia and then Rutgers. The three-year contract was extended to five years, but then I started to visit Taiwan frequently. Since Henry and I were in different places, it seemed that we could never finish the book. The publisher finally gave up on us. They did not even want us to pay back the advance. They just wanted to forget about it. When Springer learned this, it got interested in the book. We signed up with Springer, but we did not see a penny until the publication of the book in 1978. This was how we wrote our book on probability theory.

In 1962, when I was at IBM, I went to the International Mathematics Congress in Stockholm. Doob was there too. When I told him about my work with Robbins on optimal stopping, he said it was an interesting topic and suggested that we write a book for Springer. I was reluctant since we did not have enough materials—the topic was very new. However, Doob replied that being short of materials made it even more sense to write a book. Robbins agreed with Doob when I told him about the conversation, so we decided to offer a course on the subject to force ourselves to write down some notes. I offered the course, and Robbins participated too. David Siegmund was in the class. We would have a draft if we put the class notes together, but then I went to Purdue and we only met during breaks after that. It was very difficult for us to complete the book. In 1966, at the IMS Annual Meeting, I suggested to Robbins, who was then the IMS President, that since we had been lazy for two years, we should invite David Siegmund to co-author. Robbins liked the idea, Siegmund agreed and the book on optimal stopping was finally published in 1972.

**Query:** Herb Robbins was also at Purdue for a while. How did that happen?

**Chow:** The Statistics Department at Columbia was full of uncertainty when I visited there in the spring of 1965 from Purdue. I told Robbins that Haas had become a dean at Purdue and he was very interested in building up statistics. I asked if he was willing to go there, but he replied that he had to think about it. I then asked Aryeh Dvoretzky and David Siegmund, who was writing his thesis with Robbins; they both said yes. I then called Shanti Gupta. Shanti was very supportive and offered everyone a position. All of sudden, Robbins, Dvoretzky, Siegmund and I all went to Purdue. It was like a package deal.

Robbins visited Purdue for one year. Although Purdue was very interested in keeping him, Robbins

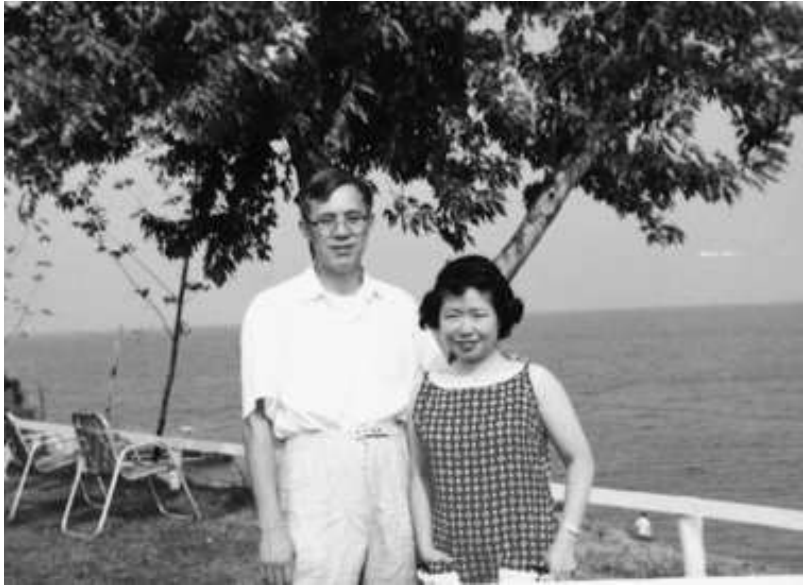

Fig. 3. *Y. S. Chow and Yi Chow in 1966.*

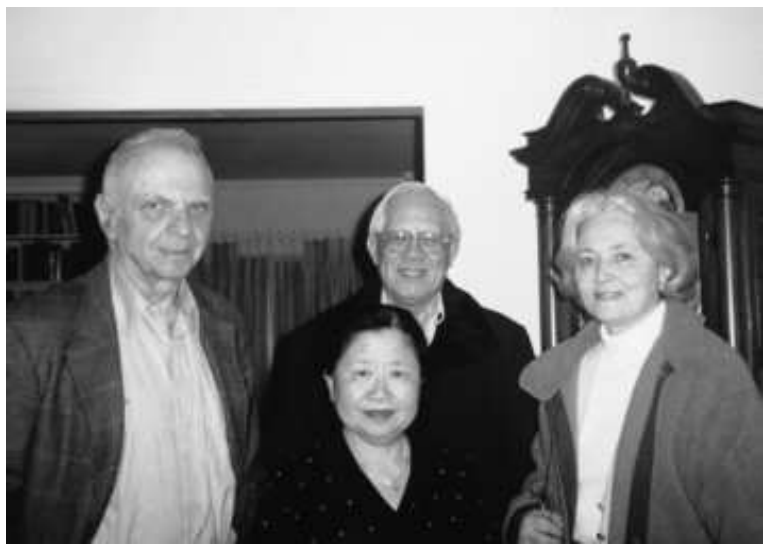

Fig. 4. *Herbert Robbins, Y. S. Chow, Yi Chow and Carol Robbins.*

did not like a school mainly composed of science and engineering. He left for the Mathematics Department at University of Michigan, which did not have a statistics department when Robbins went. Michigan wanted him to stay and formed a new department, but Robbins left anyway after helping Michigan to establish its statistics department. He visited the Electrical Engineering Department at Berkeley and then returned to Columbia in 1968. I believe that he still found New York City to be the most attractive place after his grand tour. Of course, I also came back to New York in 1968.

## ACADEMIA SINICA, TAIWAN

**Query:** You mentioned that you traveled to Taiwan often around that time. What did you do over there?

**Chow:** In 1968, Shiing-Shen Chern was recruiting scientists from overseas for Taiwan. He sent out letters asking if anyone would be willing to visit for a year. Columbia had given me credit for my service as an adjunct professor, so I qualified for a sabbatical leave in 1970, two years after my return to Columbia. I replied to Chern that I would be able to spend my sabbatical leave in Taiwan, so Chern wrote to Shichie Wang, then President of Academia Sinica in Taiwan. Again, Chern influenced my career. When Wang came to New York City, he invited my wife and me to dinner. He said he was very pleased that I would be able to go back. He would like me to serve as Director of the Institute of Mathematics for a year, as the current director had just decided to leave Academia Sinica.

My sabbatical leave turned into something far more involved and gratifying than I had anticipated. After my arrival in Taiwan, Siliang Chien became the President of Academia Sinica and he wanted to renew my appointment. I said I would consider it if he doubled the budget for the Institute. He agreed. It was a small amount, but it made the Institute of Mathematics relatively rich for a long time since the budget of individual institutes increased by a percentage every year. I served six more years as the director, from 1971 to 1976, and traveled frequently between New York and Taiwan.

In 1968, President Chiang Kai-shek consulted with physicist Luke C. L. Yuan about ways Taiwan could attract overseas Chinese scientists. Yuan suggested that the government should allow people to be in charge of scientific institutions on a part-time basis, letting them keep their permanent positions abroad, instead of demanding they choose between Taiwan and the U.S. Yuan said that such appointees would be more objective and feel free to serve to the best interest of their units. Chiang Kai-shek approved. I was one of the first to serve under this new policy.

**Query:** Could you tell us a bit about the formation of the Institute of Statistical Science in Academia Sinica?

**Chow:** My main goal as Director of the Institute of Mathematics was to promote computer science and statistics in Taiwan. In those years, Taiwan urgently needed to improve living standards. The per capita income in Taiwan was $450, and the monthly income for a professor was $100, barely enough to support a family. Professors did not have enough money to buy books, let alone to conduct research. Within the science, computer science and statistics are fields that can lead directly to greater prosperity. During my tenure at Academia Sinica, I pushed hard for these two fields. We formed research groups within the Institute of Mathematics to attract people, so that they could eventually split out to form new institutes. We tried to find ways to use Taiwanese parts to build smaller computers like the DEC PDP-8 instead of the bigger IBM type. Many well-known scientists helped our institute in this project, including some who later played important roles in the development of the computer industry in Taiwan. In 1977, shortly after I left, Academia Sinica formed the Preparatory Office for the Institute of Information Science, which formally became an institute in 1982.

From 1982 to 1987, I served as Chairman of the Board of the Preparatory Office for the Institute of Statistics, and Min-Te Chao served as Director. George C. Tiao, Gregory Chow and C. C. Li also played important roles. We invited C. R. Rao over to Taiwan so we could learn from the great success of the Indian Statistical Institute. We realized that it was not wise to copy its model exactly, since we would like to play a more active role in promoting industry and commerce. We therefore invited Moriguti from the Imperial University, Tokyo, who was instrumental in the development of the computer industry in Japan in the 1960s. He served as President of the International Statistical Institute and was very talented in administration. Moreover, I wanted the institute to maintain close ties to the government Statistical Bureau. For this, I solicited help from many people outside of statistics. The institute was formally established in 1987.

**Query:** Would you tell us something about the development of statistics in Taiwan?

**Chow:** Statistical sciences cover a wide range of areas, but we must focus on applications, on real problems. To develop statistics in Taiwan, we must rely on our own efforts and do what others are not willing to do and what others have never done. Now applied statistics is thriving in Taiwan. For example, we conduct all kinds of public polls, including election polls. Recently, Taiwan started to focus on medical science and public health, and decided to establish its National Health Research Institute, modeled after the NIH in the U.S. There are many talented biostatisticians in NIH, but they belong to different groups in various research institutes. Before the formation of Taiwan's NIH, I had a long talk with Gordon Lan. I said we must seize this opportunity to establish an independent institute of biostatistics. Taiwan's National Health Research Institute now has five divisions, including the Division of Biostatistics and Bioinformatics. Chao Agnes Hsiung is serving as Director of the division, which will eventually become an institute. It will be the first such statistical institute in the world.

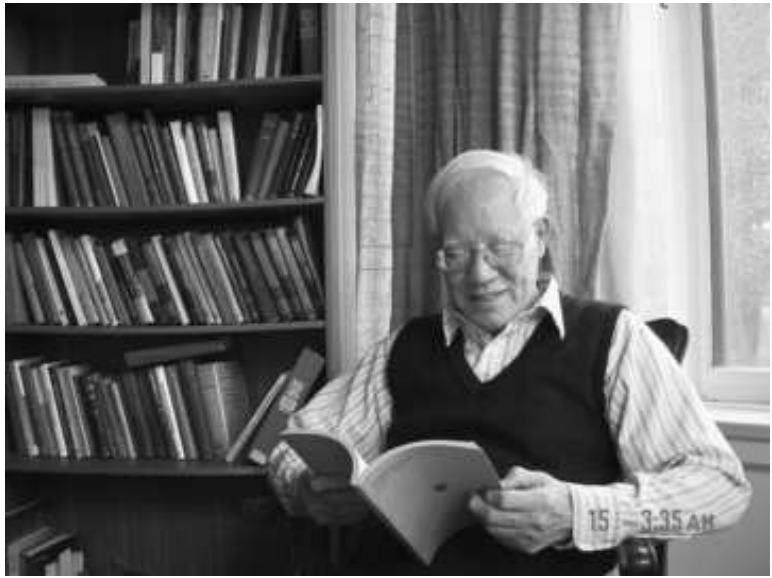

Fig. 5. *Y. S. Chow in 2003.*

We cultivated a very close relationship between the statistical community and the government in Taiwan. The Statistical Society in Taiwan is a semi-official organization, whose members are mainly from the National Statistical Bureau. The Chairman of its Board of Directors is always from the Bureau, but the majority of its directors are academics. One of the most powerful committees in the National Statistical Bureau is the Supervisory Committee on National Income, which convenes every year to approve the national income figure announced by the government. Its Chair is the General Director of Budget, Accounting and Statistics, and its members are composed of deputy secretaries of various departments (Financial Ministry, Central Bank, etc.) of the Executive Yuan (the Cabinet) and academics, in about a two-to-one ratio. The committee has to reach a unanimous decision. Some people say this is just a formality, but we have to be in the government to get to know the right people, to express our opinion and to get things done. In Taiwan, the government provides tremendous support to statistics and effective communication with it is essential. This is by no means easy to achieve, but Taiwan has been quite successful at it.

## CHINESE STUDENTS IN THE UNITED STATES

**Query:** Many Chinese students come to the U.S. to pursue their postgraduate education. You are one of the pioneers. Could you share with us your thoughts on this?

**Chow:** I would first like to share with you some background information. I belong to the first generation of students who came to the United States with financial support from American universities. Most students coming here before the revolution were supported by the Chinese government, somewhat like visiting scholars from Mainland China in the eighties. Whatever they did, whatever they said, they had to take into consideration the reaction from the government. We did not have such a close relationship with the Taiwan government, so we did not feel we had any obligation. In short, we had a different mentality.

It was extremely difficult to get scholarships from American universities, and the amount was often quite small. We had to use borrowed money to certify our financial ability when we applied for visas, and we had to find part-time jobs to support ourselves here. School authorities were well aware of this and allowed us to work on campus. The situation improved in the late fifties, partly because of the excellent academic performance of Chinese students. I remember that newspapers in Taiwan used to carry articles saying that Chinese universities had reached the world level because their graduates were directly admitted into postgraduate schools in the United States. Most students before us started from senior undergraduate courses in the U.S., including graduates from Zhejiang and Taiwan Universities, since college degrees from China were not recognized at that time.

After I joined Columbia University in 1962, I found that the Statistics Department had never admitted a Chinese student. I suggested to Robbins that Columbia offer two scholarships, one to the Mathematics Department at National Taiwan University and the other to Taiwan Normal University. They were the only two mathematics departments in Taiwan at that time, and there was no statistics department. Mathematics majors did not know that they could naturally take statistics as their specialty for postgraduate study. Robbins offered each university a four-year fellowship. Taiwan Normal sent Pi-Erh Lin, who later joined the faculty at the Florida State University, but the student nominated by National Taiwan University failed an English test at the U.S. Embassy. It turned out that National Taiwan University recommended neither the one who wanted to come most nor its best. The one who wanted to come most, Min-Te Chao, was also the best mathematics senior at National Taiwan University for that year. He spent one more year in graduate school in Taiwan before he went to Berkeley. Although he did not come to Columbia, Robbins' letter helped him to choose statistics. These two were the earliest students. From then on, many mathematics majors in Taiwan followed their path.

In 1972, Richard Nixon went to China. I saw an opportunity and suggested to Robbins that he write to Caohao Gu, my contemporary at Zhejiang University and a well-known mathematician at Fudan University. Robbins did, but because of the Cultural Revolution, Gu never received the letter. In 1979, Zhejiang University sent a delegation to the U.S. It was composed of its president and seven vice presidents, the first delegation at such a high level from a Chinese university. Since I am an alumnus of Zhejiang University, I received the delegation when it visited Columbia. We went to see Robbins after their campus tour, and Robbins offered the Mathematics Department at Zhejiang University a fellowship in the meeting upon my suggestion.

Caohao Gu was then visiting SUNY at Stony Brook. When he learned this, he asked me if the Mathematics Department at Fudan University could also have a fellowship, since it was formed by faculty from Zhejiang University. Robbins agreed to that too. Gu is a smart person—he recommended two students without ranking them, and we had no choice but to accept both. That was the year Zukang Zheng and Cun-Hui Zhang came to Columbia.

The Statistics Department at Purdue was also a pioneer in recruiting Chinese students. In 1963, when I went to Purdue, I was surprised to see eight students from Taiwan, all mathematics majors from National Taiwan University. Since I was also from there, we became good friends. Seven of them returned to Taiwan after completing their Ph.D.'s. Although I did not help them come over, I did aid their return. I take credit for the fact that the statistical community in Taiwan has a high percentage of Columbia and Purdue alumni. Sometimes I feel a little bit guilty, though, because I have remained in the States while encouraging others to return.

## ADVISING STUDENTS

**Query:** Would you share with us your experience in advising students?

**Chow:** Some people teach students how to appreciate mathematics. They describe the beauty of mathematical ideas. This is analogous to teaching swimming by asking students to watch how others swim. This is not going to work. The more you watch, the more afraid you become; you have to jump into the water. Appreciating mathematics is one thing, but actually doing it is an entirely different matter. Many influential papers are often incomplete and imperfect. For example, the proof in Birkhoff's 1932 paper on ergodic theory was not correct and the original proof of optimality of the sequential probability ratio test was also problematic. Referees nowadays may reject these papers, but they are among the best ones. It is impossible to start something with a perfect paper. We should allow our students to explore, to make a few mistakes and eventually correct them. Do not show your students too much. That will spoil them. Let them practice first and broaden their vision gradually. It is just like rolling a snowball. You have to have something small to start with. It is important to start doing research early.

We should credit students for all the good stuff in their theses and credit their advisors for the rest. Good students are made by themselves, not by their professors. A good friend of mine once said that an advisor would give students only two types of problems: those he could not solve and those he did not care to solve. I wonder why we give student problems.

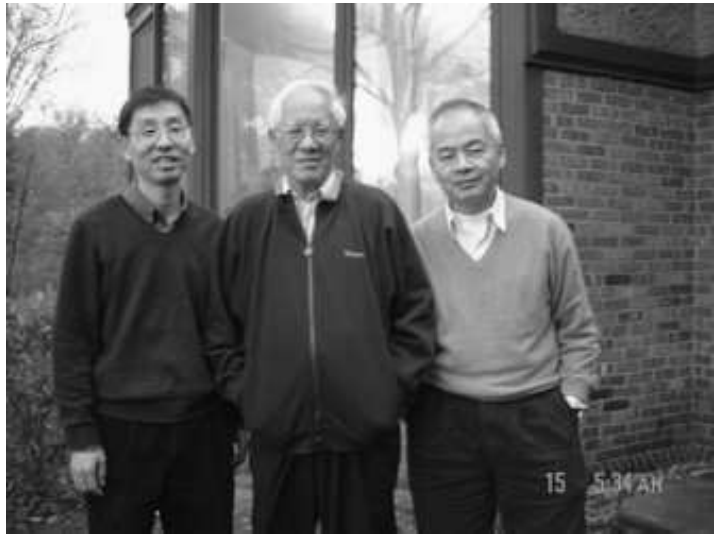

Fig. 6. *Zhiliang Ying, Y. S. Chow and Cun-Hui Zhang.*

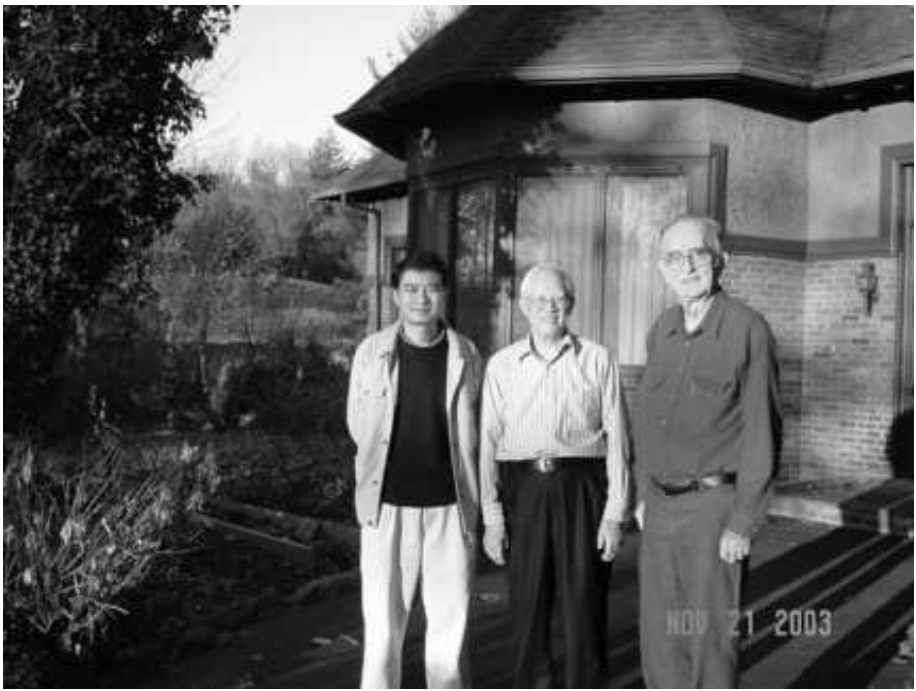

Fig. 7. *Steve Kou, Y. S. Chow and Burnett Sams.*

When I first came to Illinois, Brahana had a big influence on me. Although my college grades were not great, he said I had a strong background and forced me to take all the hard courses. I became a good student all of a sudden. This changed my life.

I realized that you could become a good student when your teacher treats you as one; if your teacher treats you as a poor student, you simply get worse and worse.

## RETIREMENT

**Query:** You retired from Columbia about ten years ago. Could you tell us something about your retirement?

**Chow:** I grew up in China and received all my early education there, for which I never paid anything. After college, I went to Taiwan. I have done something for the Chinese people in Taiwan, but not enough for those in the Mainland, to whom I owe quite a lot. I am retired, but have energy, so I have been contributing both my time and my money to education.

We, along with our friends, are subsidizing some high school students in my hometown, which is poverty-stricken. The tuition for senior high school is about the same as the average annual family income there. Some families are allotted two *mus* (about 1/3 acre) of land to farm. If their children pass the entrance exam for high school, they will have difficulties finding the money. In the past ten years, we have supported about 300 such poor students for their high school education.

We also funded two scholarship programs for college education. In the first program, we award a scholarship every year to anyone in my home county who is admitted to one of the five most selective universities in China: Peking, Tsinghua, University of Science and Technology, Fudan and Nankai. Up to now, 26 students from my hometown have been enrolled at the above universities. In the second program, we selected 20 poor students from the key high school in my hometown, and if they succeed in entering one of a group of 15 good universities, we lend a hand. Since I started the scholarship program, the number of students from this high school admitted to college has increased quite dramatically.

We also set up a math competition within the key high school in 1989, which the county government insisted on naming the "Chow Yuan Shih Math Award," despite my objections. The top five scorers in the competition each receive a prize, a small amount in the U.S., but a nontrivial figure in my hometown. Besides the prizes, I also give a banquet, so I have to go back there every year. The county magistrate, the deputy magistrate and the principal of the high school are all invited, along with the parents of the winners. The county government sends a car to pick up winners and their parents for the dinner. In addition, the local TV station always reports the results of the competition. It is a big event each year. In this way, I spend very little for very good results. In recent years, students from my hometown have always taken the top two spots in the Xiangfan Regional math competition.

Another thing in my retirement plan was to bring American-style home architecture to my hometown. Xiangfan is an ancient city, but wars stripped off all its historical landmarks. The Communists and the Nationalists fought their last battle in the city, with 25,000 casualties on each side. After the revolution, the city became desolate and required rebuilding. When I went back in 1989, I saw houses along a street 5 kilometers long, all of the same height (about 7 meters), all of the same grey color and none with a bathroom. I could not believe my eyes! Those who need the toilet have to go down the street to the public latrine. It was impossible to explain to the locals what could be done, so I decided to build a house myself. I modeled my house after Illinois' Arlington Park and brought architects to modify my design. I built the house using Chinese products no matter how much I had to spend, and hired local builders no matter how inefficient they were. Nevertheless, the house was finished in three years. It is beautiful, but it needs repairing often. When I pay for replacement materials, I complain loudly and always ask merchants to stock better stuff. Now, when the locals build their own houses, they include bathrooms too. What's more, I used a light blue color, and it has become a fashionable color for houses in my hometown. I also used parquet hardwood floors, which was new to the locals, but has become very popular now.